\documentclass[11pt]{article}
\usepackage[utf8]{inputenc}
\usepackage{amssymb}
\usepackage{amsthm}
\usepackage{graphicx}
\usepackage[colorlinks=true,citecolor=blue,urlcolor=blue]{hyperref}

\def\endproof{\vrule height6pt width6pt depth0pt}

\theoremstyle{plain}
\newtheorem{thm}{Theorem}[section]
\newtheorem{lem}[thm]{Lemma}
\newtheorem{prop}[thm]{Proposition}
\newtheorem*{cor}{Corollary}

\theoremstyle{definition}
\newtheorem{defn}{Definition}[section]
\newtheorem{conj}{Conjecture}[section]

\theoremstyle{remark}

\pdfoutput=1

\begin{document}
\title{Orthogonal Representation of Graphs\\
(Working paper)}
\author{Alberto Solís Encina \\Jose Ramón Portillo
\\\vspace{6pt} Departamento de Matemática Aplicada I, \\Universidad de Sevilla, España}
\date{\today v. 1.2}
\maketitle

\begin{abstract}
Orthogonal Graph Representations are essential tools for testing existence of hidden variables in quantum theory. As required by the interpretation of Copenhaghe on the foundations of quantum mechanics, a physical observable is not determined before its observation. Conducting experiments \emph {quantum contextuality} or the information capacity of a quantum system are closely related to the orthogonal representations.
\end{abstract}

\section{Definition of Ortogonal Representation (OR)}
\begin{defn}
Let $ G = (V, E) $, an orthogonal representation (hereinafter OR) is a function $ \sigma: V \rightarrow \mathbb {R} ^ d $, with $ d $ the dimension, so that checks the next:
\begin{itemize} 
\item If two vertices $ i, j $ are adjacent $ \rightarrow \langle \sigma (i), \sigma (j) \rangle = $ 0. 
\item $ \sigma (i) \not = \sigma (j) $.
\end{itemize}
This representation is called \textit{orthonormal} if verified  $ || \sigma (v)|| = 1 $ for all $ v \in V (G) $ and  is \textit{minimum} if there is no representation for $ G $ with $ d '<d $. At this minimum dimension $ d $ is known as \textbf{orthogonal range} of the graph. $ \rho (G) $.
\end{defn}
The assignment $ \sigma \equiv $ 0 is a trivial orthogonal representation for any graph, the same way that for any graph $ G $ such that $ | V (G) | = n $ always has a representation (does not have to be minimal ) $ \sigma (G) \in \mathbb {R} ^ n $. The main difficulty of the problem is the minimization of the dimension because the space of the representation in vector  should be optimizated and also the dimension is not unique. 
\paragraph{}
The problem of orthogonal representation, from the purely mathematical point of view, was approached by Lovasz ~\cite{1} where he established a necessary and sufficient condition for a specific type of representation known as \emph{general position}. Subsequently many authors ~\cite{12,13,14} described other representations and bounds for the minimum dimensions of them. However, in the context of studies on quantum bases discussed above, another representation appears ~\cite{3,7,10,11} for graphs: \emph{faithful} representations (hereinafter FOR).
\begin{defn}
Let $ G = (V, E) $, a faithful orthogonal representation (hereinafter FOR) is a function $ \rho ^ \perp: V \rightarrow \mathbb {R} ^ d $, with $ d $ the dimension such that verifies the following:
\begin{itemize} 
\item If two vertices $ i, j $ are adjacent $ \leftrightarrow \langle \sigma ^ \perp (i), \sigma ^ \perp (j) \rangle = $ 0. 
\item $ \sigma ^ \perp (i) \not = \sigma ^ \perp (j) $.
\end{itemize}
Similarly, this representation is called \textit {orthonormal} if verified  $ || \sigma ^ \perp (v)|| = 1  $ for all $ v \in V (G) $ and is \textit {minimum} if there is no representation for $ G $ with $ d '<d $. At this minimum dimension $ d $ is called {\bf orthogonal range faithful of graph:} $ \rho ^ \perp (G) $.
\end{defn}
The inclusion of the necessary and sufficient condition restricts more the problem of representation, so much so that Lovasz ~\cite{1} only gives a characterization of the dimension for his own orthogonal representation.
In the main theorem in ~\cite{1} Lovasz defines a characterization for a type of representation that he called \textit{orthogonal Representation in general position} (hereinafter ORGP).
\begin{defn} Let $ G = (V, E) $, one ORGP is a function $ \sigma_ {pg}: V \rightarrow \mathbb {R} ^ d $, with $ d $ the dimension, so which verifies the following:
\begin{itemize} 
\item If two vertices $ i, j $ are adjacent $ \rightarrow \langle \sigma_{pg}(i) \sigma_{pg}(j) \rangle =  0$. 
\item Every vector set of size $ d $ is linearly independent. 
\item $ \sigma_ {pg} (i) \not = \sigma_ {pg} (j) $.
\end{itemize}
Similarly, this representation is called \textit {orthonormal} if verified  $ || \sigma_ {pg} (v)|| = 1 $ for all $ v \in V (G) $ and is \textit {minimum} if there is no representation for $ G $ with $ d '<d $. At this minimum dimension $ d $ is called {\bf orthogonal range in general position of the graph:} $ \rho_ {pg} (G) $.
\end{defn}
We have omitted from the definition that Lovasz considers orthogonal vertices those which are not adjacent. But the result is the same using the complementary graph $ \bar {G} $.
\paragraph{}
The relationship of the parameters studied by Lovasz in their work are connected with quantum theory. Adan Cabello et al. ~\cite{10} show that there are experiments of quantum measuring contextuality (determined by their inequality equations) which can be represented as graphs and also the number thereof Lovasz sets a limit for the correlation value given by quantum theory in these experiments. \\
The Importance of $ \theta (G) $ ~\cite{8} defines a new type of essential orthogonal representation in quantum mechanics.
\begin{defn} Let $ G = (V, E) $, one Lovasz-optimal faithful orthogonal representation (hereinafter FORLO) is a function $ \sigma_ {pg}: V \rightarrow \mathbb {R} ^ d $, $ d $ with the dimension, so that verifies the following:
\begin{itemize} \item If two vertices $ i, j $ are adjacent $ \rightarrow \langle \sigma_ {so} (i) \sigma_ {so} (j) \rangle = $ 0. 
\item The representation obtained get the Lovasz number $\theta$. 
\item $ \sigma_{pg}(i) \not = \sigma_{lo}(j) $.
\end{itemize}
Similarly, this representation is called \textit {orthonormal} if verified  $ || \sigma_ {so} (v)|| = 1 $ for all $ v \in V (G) $ and is \textit {minimum} if there is no representation for $ G $ with $ d '<d $. This minimum dimension $ d $ is called \textbf{Lovasz-optimal orthogonal range of graph: } $ \rho ^ \perp_{lo} (G) $.
\end{defn}

\section{OR from the geometrical point of view}
Given a representation with base assigned to a set of vertices V = $\left\lbrace  v_1, ..., v_i, ..., v_d \right\rbrace $ in dimension $d$ forming the structure of a complete graph $G = K_d $, get a FOR induced by subgraph $ G' \subset G $ is done by removing edges. These operations correspond geometrically to perform rotations about the set of vector $V$, so that couples $ (v_i, v_j) \in \bar {E} (G) $, with $ i \not = j $, are selected and rotate one of the vector around a given axis.
In each step of these procedure is logical that the restrictions imposed by $ G'$ can not be representable for dimension $ d $, so the dimension may fluctuate to a greater or lesser in every step.
\begin{prop} Given any graph $ G $. If $ G $ is representable in $ \mathbb{R}^d $, then there exists a transformation in the space $ T $ on the canonical basis $ B (\mathbb{R}^{| V (G) |}) $ of vectors representing $ v_i \in V(G) $, so that $T = G_{v_1}\cdot G^1_{v_i} \cdots G^k{v_{|V(G)|}}$ is a sequence of turns on each vector and $T(B(\mathbb{R}^{|V(G)|})$ is an FOR of $ G $ in dimension $d$.
\end{prop}

{\em proof:}
At every turn does not have to keep the dimension $ d $, in fact may increase or decrease. The reason for decrease in one dimension to make one round (remove an edge in the graph) or should instead increase is due to the \textit{distinction between vectors}, since the vectors must be different for different vertices because proportionality are not allowed and therefore there are cases
in which a transformation to eliminate a component in the vector representation (the same for all), but in other cases the corresponding transformation never keeps the condition of distinguishability of vertices and therefore appears necessary to add a new component to the group. Assuming a representation, there is a transformation $ T $, do not know if the rotation are commutative, which gets the FOR.
$$
T \left( \begin{array}{ccc}
v_{11} & \cdots & v_{n1} \\
\vdots & \ddots & \vdots  \\
v_{1n} & \cdots & v_{nn} \end{array} \right) = 
\left( \begin{array}{ccc}
v'_{11} & \cdots & v'_{n1} \\
\vdots & \ddots & \vdots  \\
v'_{1d} & \cdots & v'_{nd} \end{array} \right) = V'
$$
Clearly, the rotation is unique because it supports a rotation of $ x + 2k \pi $.
\endproof
\section{Banned Graphs} 
Banned graphs are graphs with the particular property that monopolize the space $ \rho ^ {\perp} = n $ without becoming a $ K_n $, is due to restrictions on the conditions of FOR. \\
Some previous definitions:
\begin{defn} A graph $ G $ is said \textbf{critical} for  dimension $ d $ when the dimension change to remove any edge in $ G $.
\end{defn}
\begin{defn} A \textbf{Banned Graph} $ G_{fb} $, is the graph critical for dimension $d$ such that $\rho^\perp (G_{fb}) = d$. It said to be banned for dimension $d-1$.
\end{defn}
Banned graphs provide the best known lower bound for $ \rho ^ \perp $, with an overall dimension for all graphs:
\begin{lem} Be a graph $ G $ and $ G_{fb} $ the largest graph induced prohibited. \\  
$$\rho^\perp(G_{fb}) \le \rho^\perp(G) \le |V(G)|$$
\end{lem}
The existence of Banned graphs comes from studying the situation in which we try to force a dimension in a graph.
Mainly due to the need to maintain 2 properties in FOR:
\begin{enumerate} 
\item \textbf{Uniqueness}. That is, no duplication vectors allowed. 
\item \textbf{Orthogonality}. The graph should reproduce, in a  necessary and sufficient way, the set of resultant vectors.
\end{enumerate}
These two properties of the FOR define two different types of structures that can be found prohibited.
\subsection{Banned Graph by duplicity}
When in a graph are 2 vertices $ v_i $ and $ v_j $ (with $ i \neq j $) that perfrom the same function, ie, they have the same neighbors then they are, inside the structure, indistinguishable. This situation, when
perfrom a projection on the space assigning vectors to the vertices FOR able to cause $ v_i \| v_j $. These structure is one of the subfamilies in family type I of banned graphs: $ F_n = K_n - \lfloor n / 2 \rfloor $, which actually consists of removing a maximal matching of $ K_n $: $ K_n - M (K_n) $. So there are always at least 2 vertices with duplicity. The particularity of this subfamily is that $|V(F_n)| = \rho^perp(F_n) = n$ \\
Examples are shown in the following figures:

\begin{figure} [htbl!]
\begin{minipage} [b] {0.5 \linewidth}
\centering 
\includegraphics[width = 0.3 \textwidth]{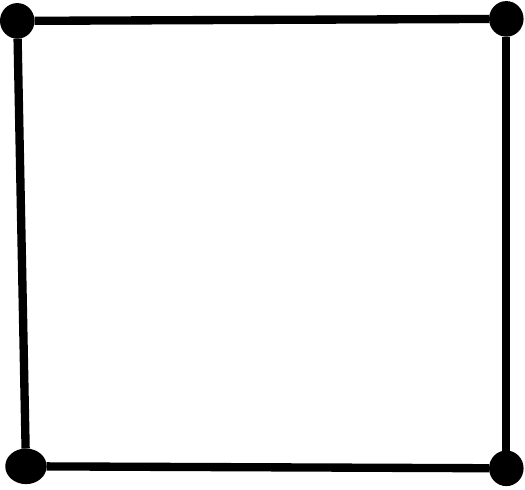} \caption{Graph not realizable in dimension 3} \label{label}
\end{minipage}
\hspace{0.5cm}
\begin{minipage} [b] {0.5 \linewidth}
\centering 
\includegraphics[width = 0.3 \textwidth]{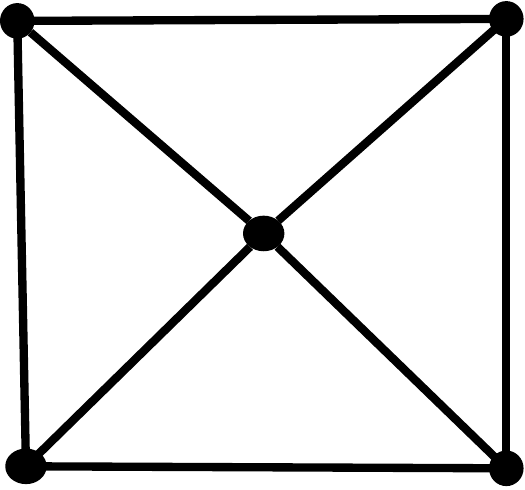} \caption{Graph not realizable in dimension 4} \label{label}
\end{minipage}
\end{figure}
\begin{figure} [htbl!]
\begin{minipage} [b] {0.5 \linewidth}
\centering 
\includegraphics[width = 0.5 \textwidth]{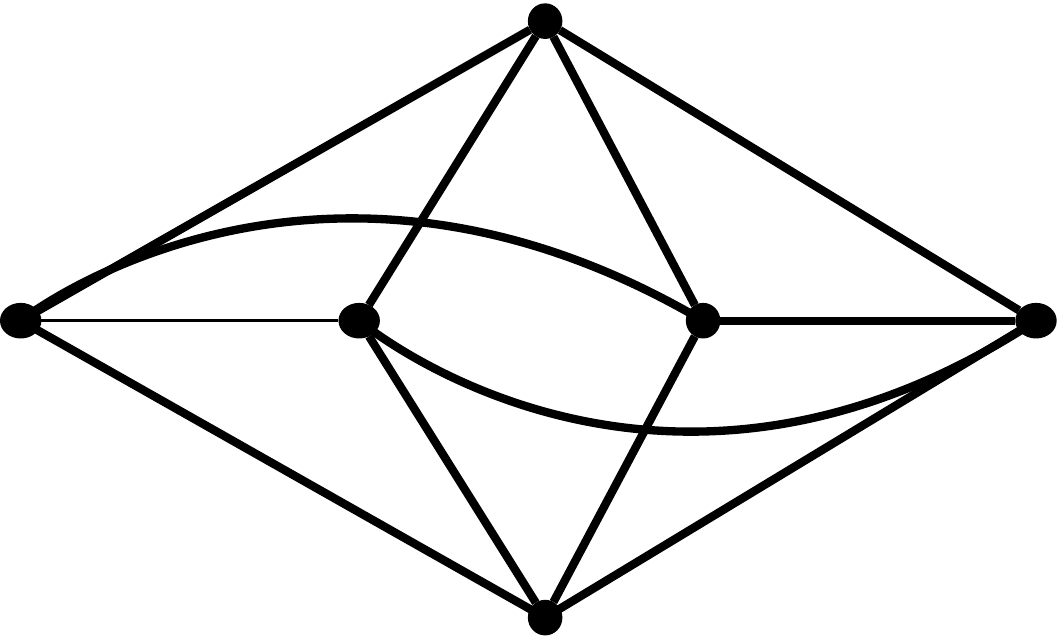} \caption{Graph not realizable in dimension 5} \label{label}
\end{minipage}
\hspace{0.5cm}
\begin{minipage} [b] {0.5 \linewidth}
\centering 
\includegraphics[width = 0.5 \textwidth]{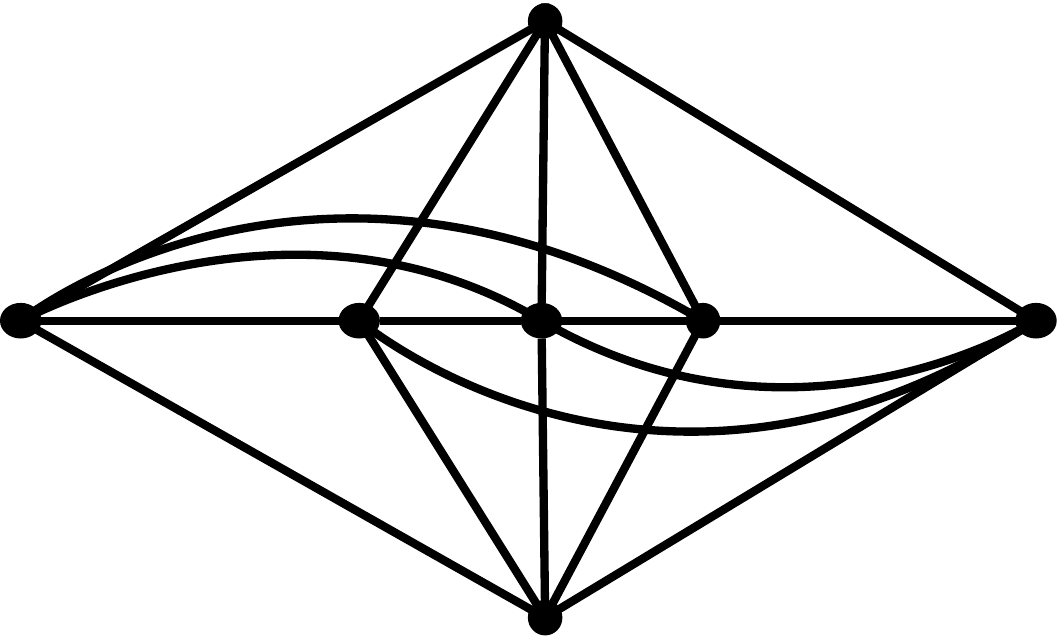} \caption{Graph not realizable in dimension 6} \label{label}
\end{minipage}
\end{figure}
\begin{conj} If we have a banned graph  to dim $ k $ and doubled a vertex of valence $ k $, then I have a graph banned in dim $ k + 1$.
\end{conj}
\textcolor {red} {No demo ...}
\paragraph{Properties of banned graphs type I with $ | V (G) | = d $}
Let $ F_d $ the banned graph  representable in dimension $ d $ and banned for$ d- 1$.
\begin{itemize} 
\item $ n $ (command) $ = d $. 
\item $ e $ (size) $ = \lceil d ^ 2/2-d \rceil $. Follow the \textit{sequence} $ Eliptic troublemaker R_n (2,4) $. 
\item Eccentricity: $ \xi (v), \forall v \in V (F_d) \le $ 2. 
\item Diameter: $ max \{\xi (v) \} = $ 2. 
\item Girth (length of the shortest cycle): $ 4 $ if $ n \ge $ 4. In another case $ \infty $. 
\item Radio: $ rad (G) = min \{\xi (v) \} \in [1,2] $, depending on whether it is even or odd.
\end{itemize}

We get the simplest structure of type I forcing two vertices to be orthogonal without direct adjacency, that is, using a structure that occupies space such that both vertices have the same \textit{behavior} into the graph.
\begin{figure} [h!]
\centering
\centerline {\includegraphics[scale = 0.40]{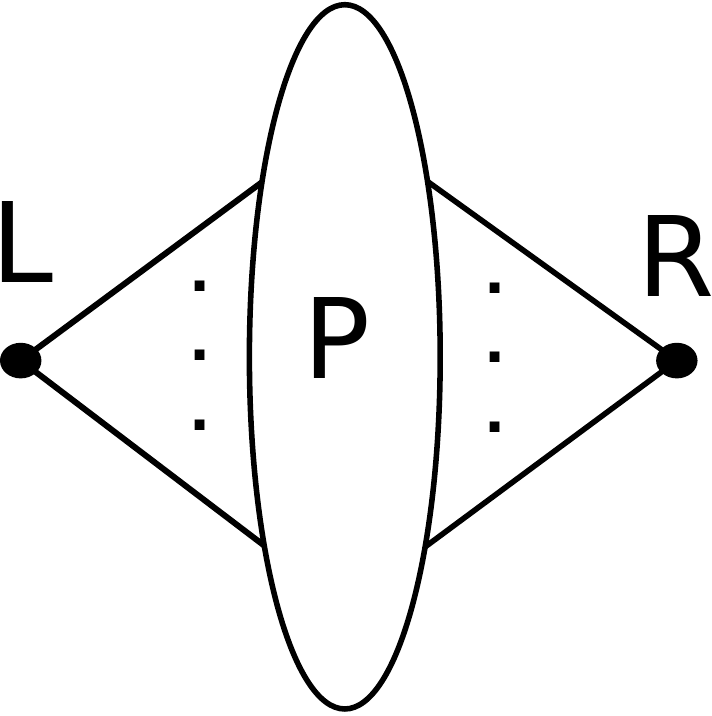}}
\caption {\label{Figa}
The structure of hidden edge in dim $ d $. P is a graph such that $ \rho ^ \perp (P) = d $. $ G = P \cup {l} \cup {r} $ turns out to be a graph, if we try to continue to represent dim $ d $ or $ d + 1 $, only got $ l \bot $ r. 
}
\end{figure}
\\
If in the above figure we consider that $ P $ is the graph banned to dim $ d $, then $ G = P \cup {l} \cup {r} $ be the graph is banned to dim $ d + 1$.
$$\rho ^\perp(G_{typeI}) = \rho ^\perp (P) +2$$
\paragraph{}
There are other graphs of type I different of above structures but are also critical.
This family is similar to Figure 6, but the are not complete adjacency between the vertices.
\begin{figure} [h!]
\centering
\centerline {\includegraphics[scale = 0.40]{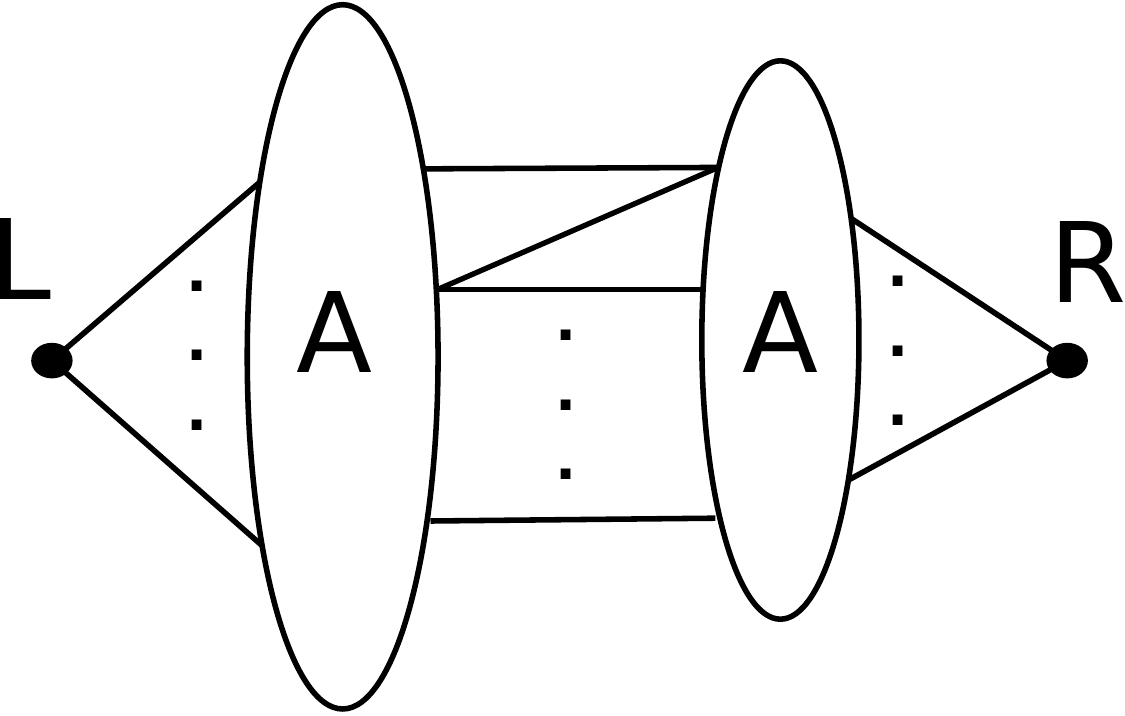}}
\caption {\label{Fig6}
Another structure of Type I. $A$ is a duplicate graph and the are not complete adjacency between the vertices $l$ with $A$ and $r$ with $A$.
}
\end{figure}
\\
To explain this case we know a family of graphs based on multiple cycles of length 4: $ C_{4n} $, with $ n \ge $ 2, and does not contain as banned subgraph any of the previously studied.
\begin{figure} [htbl!]
\begin{minipage} [b] {0.5 \linewidth}
\centering 
\includegraphics[width = 0.5 \textwidth]{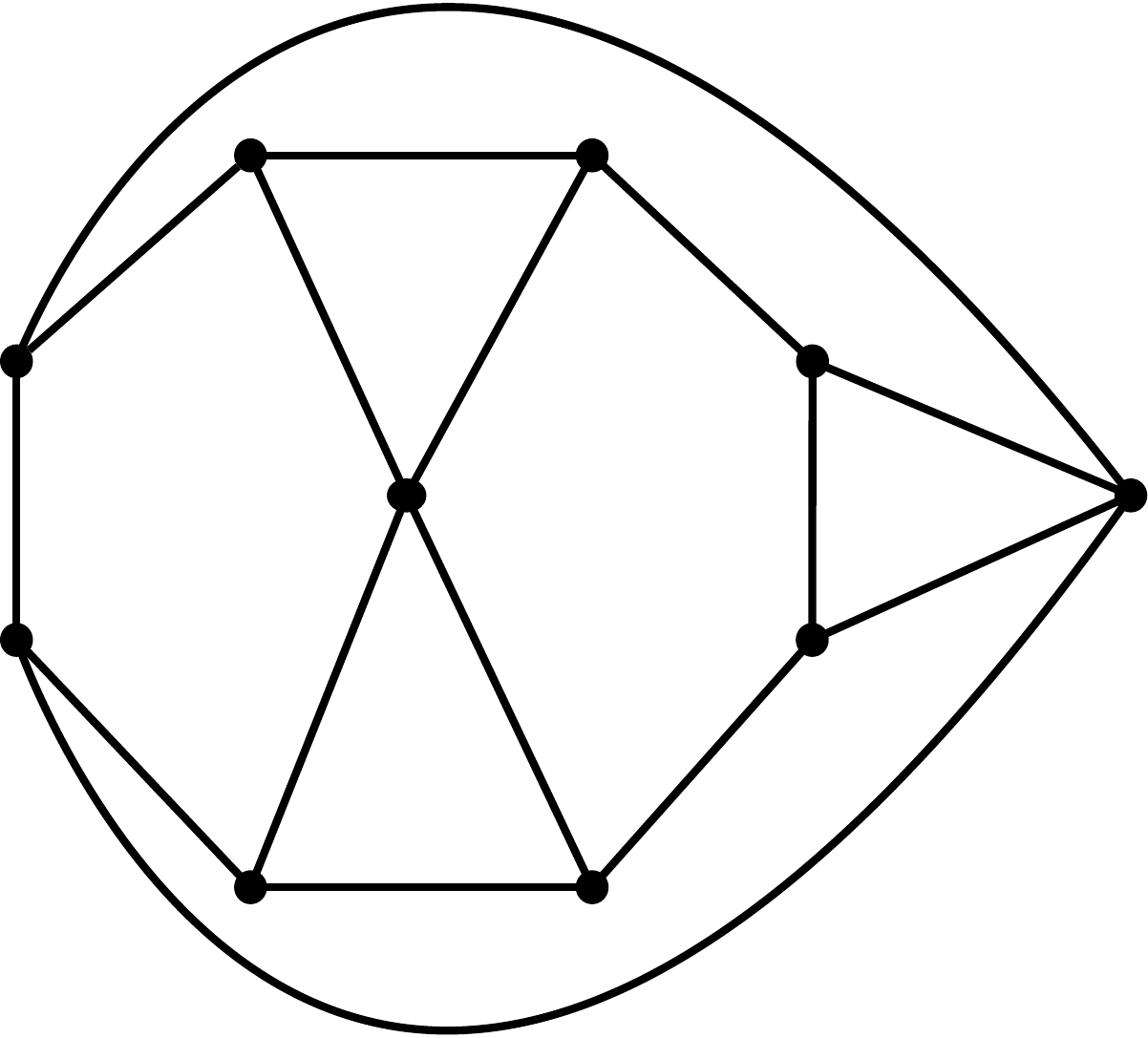} \caption{Graph with $ C_8 $. This graph is not realizable in dimension 3} \label{label}
\end{minipage}
\hspace{0.5cm}
\begin{minipage} [b] {0.5 \linewidth}
\centering 
\includegraphics[width = 0.5 \textwidth]{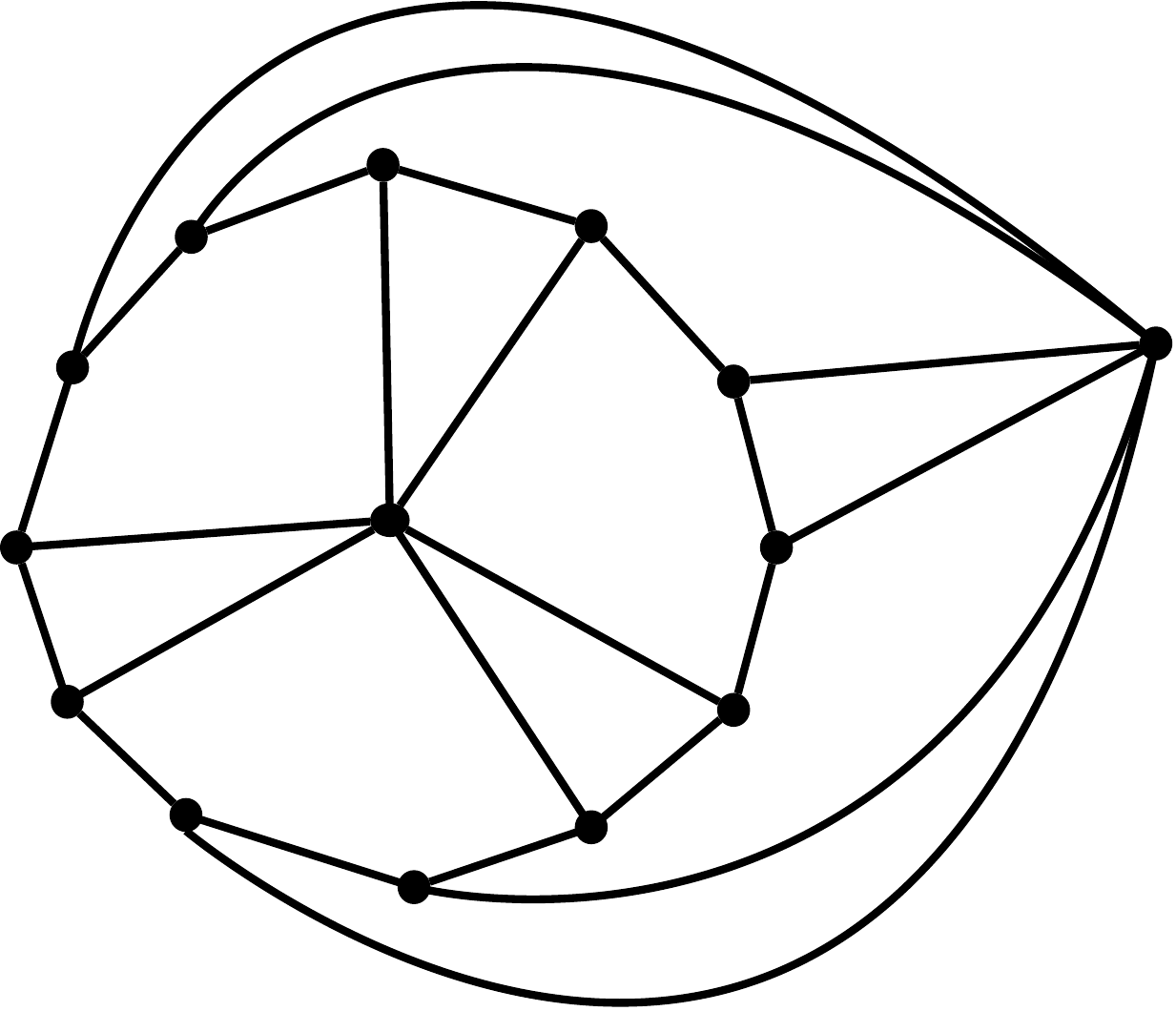} \caption{Graph with $ C_ {12} $. This graph is not realizable in dimension 3} \label{label}
\end{minipage}
\end{figure}

\subsection{Banned Graph by edge hidden}
\begin{defn} A \textbf {hidden} edge is a virtual orthogonality result of the difference between a graph $ G $ and a vector representation of its vertices. It occurs when two vertices $ v_i, v_j $ verify: 
\begin{itemize} 
\item $ (v_i, v_j) \not \in E (G) $ 
\item $ v_i \bot v_j $ 
\end{itemize}
\end{defn}

This is the structure of type II with orthogonality by cross assignment.
\begin{figure} [h!]
\centering
\centerline{\includegraphics[scale = 0.40]{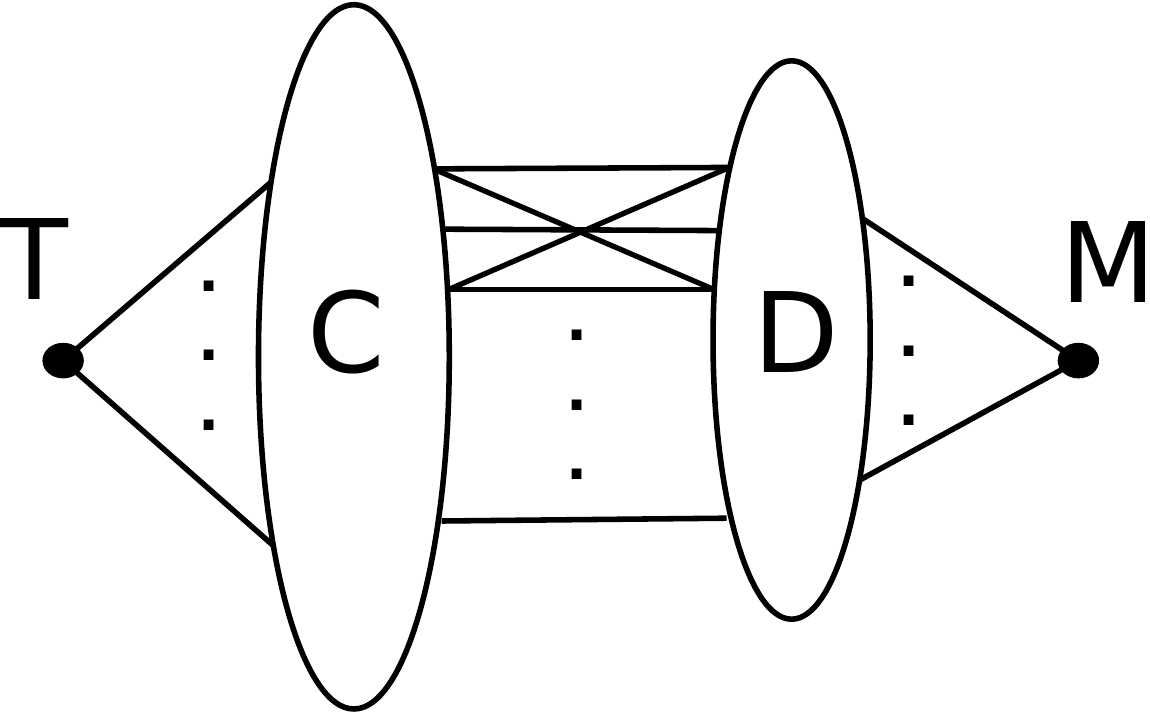}}
\caption {\label{Figh}
Structure type II.  Hidden edge in dim $ d $ with dual core.
}
\end{figure}
The diagram in Figure 9 with dim $ d $ is: $ D $ is a linear subspace of dimension $ k $, $ C $ is another linear manifold of dimension $ d-k $ and sets $ T $ and $ M $ are linear subspace each containing only one vertex $ t $ and $ m $, respectively. Note that $ T \bot C $, as $ t \in D $ this means that $ T \bot M $, hence $ t $  and $ m $ are orthogonal and there is a hidden edge between them. If this happens we know it is not possible to assign vectors in dimension $ d $ in this graph.
\\
The underlying recursive feature in the Banned graphs is well noted. Scale not only adding vertices in accordance with the above structure, but the total binding of two graphs also generates another banned graph.

$$\rho ^\perp(G_{typeII}) = \rho^\perp(C \cup D) + 1$$

\begin{prop} Be a graph $ G $ with structure hidden edges as in Figure 9 in dimension $ d $, with $ C, D \in G $. Then $ C \cup D $ is the graph banned type I for dimension $ d-1 $ more an edge when $ | C | $ and $ | D | $ are odd. Otherwise prohibited graph matches for dimension $ d-1 $.

\end{prop}
{\em proof:}
It is sufficient to calculate the difference between the edges of the graph generated according to the scheme of Figure 6 and those belonging to the banned graphs. For the case where $ C $ and $ D $ has an odd number of vertices is 1 the difference edge. For all other cases the resulting graph is the banned graph for dimension $ d-1 $.
\endproof
\paragraph{}
These families are discussed above be banned graphs with the minimum number of vertices, ie, the smallest graph is  represent necessary $ d $ dimension is the dimension graph prohibited for $ d-1 $. \\
\begin{conj} Any graph Banned critical to dim $ d $ has at least 2 vertices \textit {'isomorphic'}. NOTE: 2 isomorphic vertices are those that removing either (either) the graph has smaller dimension (being critical) and also both generated subgraphs are isomorphic.
\end{conj}
\textcolor{red}{No demo ...}
\begin{conj}
If $ G $ is a critical graph (not banned) and $ \rho ^ \perp (G) = d $, then $ G - M (G) $ (with $ M $ a maximal match) has $ \rho ^ \perp (G - M (G)) \le d-1 $ 
\end{conj}
\textcolor{red}{No demo ...}
\paragraph{}

\paragraph{Banned Graphs Family}

Mathias proposed a group of Banned graphs. We handle categorize by type (type I or II). However, is important to think about this clasification first of all.
First, every graph has a banned dimension, ie $ \rho ^\perp (G) = \rho ^\perp (G)  -1$, clearly. This pre-orthogonal dimension range of the graph is insufficient to represent vectors, the reason may be based on arguments of type I or type II. It is like saying that Petersen is Banned to dim 2, in fact it is because $ \rho ^\perp (P) =  3$. \\
Perhaps most important-difficult-essential thing is to identify when a dimension is critical to a graph, by this I mean that although we change the dimension that we tested, the morphology of the graph remains the same. Therefore, elements such as duplication, that is, that two vertices are "indistinguishable" (sharing neighbors) is a property that does not vary with the size, although it may itself depend on the number of adjacencies.
Therefore we clearly identify a graph with type I or type II is to assert that there is or not parallel (type I) or hidden edge (type II) in the graph. But this requires existence graph quantifying parameters to determine the existence of the above elements on the specific graph.
What is special about the family we know? The most important family in our paper is type I and also has that $ | V | = \rho ^\perp (G) $.
Owned by the dimension of a graph we used to obtain higher levels of graphs. However, what happens when $ | V |> \rho ^\perp (G) $? Because it is possible found other families of banned graph with the same structure of type I.
\paragraph{}
\textcolor{red}{Open question:} Be a graph $ G = (V, E) $ and dimension $ d $, Is G type I for $ d $?
Suffice to note that, for $ d $ given, there are two vectors with the same representation. But we do not know how to determine the $ d $
\paragraph{}
You can create a table by specifying the family of graphs banned in the number of vertices and the banned dimension.
$$\mathbf{F_B} = \{\rho^\perp_{d},\rho^\perp_{d+1}, \rho^\perp_{d+2}, \cdots, \rho^\perp_{\infty}\}$$
The table would be:
\begin{center}
\begin{tabular}{c|c|c|c|c|c|c}
& $\rho^\perp_{d+1}$ & $\rho^\perp_{d+2}$ & $\cdots$ &$\rho^\perp_{d+7}$ & $\cdots$&$\rho^\perp_{\infty}$ \\ \hline
$d=3$ & $C_4$ & $\not \exists$ & $\cdots$ & $C_5 + C_3$& $\cdots$&? \\ \hline
$d=4$ & $R_5$ & ? & $\cdots$ & $C_5 + C_4$& $\cdots$&? \\ \hline
$d=5$ & ... & ? & $\cdots$ & ? & $\cdots$&? \\ \hline
\end{tabular}
\end{center}

The $ \rho ^ \perp_ { d + 1}$ family is the banned graphs $K_n - \lfloor n/2 \rfloor$.
\section{Known dimensions} 
\begin{itemize} 
\item $ \rho ^ {\perp} (K_n) = n $.
 \item $ \rho ^ {\perp} (F_n) = n $ and $ F_n = K_n - \lfloor n / 2 \rfloor $ hidden edges. 
 \item $ \rho ^ {\perp} (C_n) = 3 $, unless $ \rho ^ {\perp} (C_4) = $ 4. 
 \item two-tree: $ \rho ^ {\perp} (T_2) = $ 4. Two dimension using $ \rho ^ {\perp} (G) \ge mr ^ {+} (G) $ ~\cite{41}. 
 \item $ \rho ^ {\perp} (Petersen) = $ 3. $$ FOR (Petersen) = \left (\begin {array} {cccccccccc}
-1 & 1 & 1 & 1 & 0 & -1 & -1 & 0 & -1 & -1 \\
-1 & -1 & 1 & -1 & -1 & 0 & -1 & -1 & 1 & 0 \\
-1 & 0 & 0 & -1 & 1 & 1 & 1 & -1 & -1 & -1 \end {array} \right) $$ 
\item $ \rho ^ {\perp} (J (5,2) = \bar {P}) \in [5,7 ] $. 
\item $ \rho ^ {\perp} (K_ {n, m}) = $ 4. By having $ C_4 $: $ 4 \le \rho ^ {\perp} (K_ {n, m}) \le n + m $. It can be shown that there is always a representation of $ K_ {n, m} $ in dim 4 using two orthogonal sets are together where all elements are \textit {'coplanar'}. 
\item $ \rho ^ {\perp} (K_ {n, n} - e) = 5$. Structure type II.

\item $ \rho ^ {\perp} (\bar {C_n}) = n-2 $. Josera showed that $ \rho ^ {\perp} (\bar {C_n} \le n-2 $, with n odd. We know from ~\cite{43} that $ \rho ^ {\perp} (G) \ge  mr ^ {+} (G) \ge mr (G) $ As $ mr (\bar {C_n}) = n-2 $ is shown 
\item Paley Graphs. (proven computationally up to 29) $ P ( q) $, with $ q $ prime, verifies that all orthogonal vectors intersect at the same angle and $ \rho^{\perp} (G) = \lceil q / 2 \rceil $. \\For $ P (q ^ 2) $ only have the example of 9 and 25, although we believe that $ 25 = 5 ^ 2$ out 3 angles. 
\item \textcolor {red} {No demo ...} We know that $ \rho ^ \perp \le \omega (G) + $ 1 (because graphs do not contain prohibited maximum clique but simulate As $ \omega (G) \le | G | -. \alpha (G) + $ 2 ~\cite{45} $$ \rho ^ \perp (G) \le | G | - \alpha (G) + 3 $$ 
\item For $ G $ vertex-transitive self-complementary and $ \rho ^ \perp (G) = \lceil n / 2 \rceil $. 
\end{itemize}
\paragraph{}
\begin{thm} Let $ G $ a graph anyone, $ VC $ the set of vertices cutting $ \bar{G} $, $ n = | V (G) | $ and $ m = | VC | $. 
$$\rho^\perp(\bar{G}-VC) \le \rho^\perp(\bar{G}) \le \rho^\perp(\bar{G}-VC)+m \le n-m$$ 
\end{thm} {\em proof:} Trivial using Lovasz dimension. \endproof

\paragraph{Regular graphs:} We know $ 2 \delta (G) + 2-n \le K (G) \le \Delta (G) $ ~\cite{42}. Using that if a graph $ G $ is regular, then $ \bar {G} $ is n-r-1 regular. 
$$ 2 (n-r-1) + 2-n \le K (\bar {G}) \le n-r-1 $$ 
$$ n-2r \le K (\bar{G}) \le n-r-1 $$ 
theorem by Lovasz, k-connected $ = n-d $, ie, $ d  = n-k$. Substituting: 
$$ n-2r \le n-d \le n-r-1 \rightarrow  2r \ge d \ge r + 1 $$ 
This result holds for representation in general position. Dimensional with FOR would be: $$ \rho ^ {\perp} (G) \le \rho_ {pg} (G) \le 2r $$ 

\paragraph{Self-complementary Graph:} Let $G=(V,E)$ self-complementary and $|V| = n$, we know
$$\rho ^\perp(G) + \rho ^\perp(\bar{G}) = 2 \rho ^\perp(G) \geq n$$
\begin{prop}
Let $G=(V,E)$ self-complementary and $|V| = n$, then
$$\rho ^ \perp(G) \geq \lceil n/2 \rceil$$
\end{prop} {\em proof:} Let $K_n$ and 
 $
\sigma ^\perp (K_n) =  \left( \begin{array}{ccc}
e_{11} & \cdots & e_{n1} \\
\vdots & \ddots & \vdots  \\
e_{1n} & \cdots & e_{nn} \end{array} \right) 
$
\\
Let $G$ self-complementary subgraph of $K_n$. $G$ have $n(n-1)/4$ edges, ie, the half restrictions of $K_n$ therefore lost $n(n-1)/4$ ortogonalities. In the best case for remove that ortogonalities we can remove up to $n/2$ components in representation $\sigma^\perp(K_n)$. $\sigma^\perp(G)$ is obtained as subset (submatrix) of $\sigma^\perp(K_n)$
 \endproof 

\begin{conj}
If $G=(V,E)$ is self-complementary and vertice-transitive, then
$$ \rho ^\perp(G) = \lceil n/2 \rceil$$
\end{conj}
{\em proof:} We know:
\begin{itemize}
\item If $G$ is self-complementary and not vertice-transitive may be false.
\item If $G$ is vertice-transitive and not self-complementary there is counterexample: $\rho^\perp(C_{18})=4$ and $\frac{|V(C_{18})|}{2} = 9$. 
\end{itemize}
We think that in this cases we have the best case in above proposition, such that the group
 of $n/2$ vectors which remove ortogonalities have redundancy in the graph and we can use vectors by linal combination of the rest of the graph, so we reduce the dimension of the graph to $n/2$.
 \endproof 
 
 \begin{cor}
 Let $G$ self-complementary
 $$\rho ^\perp(G) + \rho ^\perp(\bar{G}) \geq |V(G)|$$
 \end{cor}
 
\paragraph{Distance between $\rho^\perp(G)$ and $\rho ^\perp(\bar{G})$:} We know 
$$K(\bar{G_c}) \leq K(\bar{G}) \leq K(\bar{G-C} + K(G)) \equiv min \left\lbrace n(G-C-H_i)\right\rbrace + K(G)$$
by ~\cite{44}. If $G$ have dimension $d$, ie, $\rho ^\perp_{pg}(G)=d \Leftrightarrow K(G) = n-d$. Let $\bar{G}$ where $\rho ^\perp_{pg}(\bar{G})=d' \Leftrightarrow K(\bar{G}) = n-d'$. Then
$$d' \geq d - min \left\lbrace n(G-C-H_i)\right\rbrace  \Rightarrow d-d' \leq min\left\lbrace n(G-C-H_i)\right\rbrace  \leq n(G-C)$$
\begin{prop}
If $\rho^\perp(G) > \rho^\perp(\bar{G})$ then
$$\rho^\perp(G) - \rho^\perp(\bar{G}) \leq n(G-C)$$
 \end{prop}
 
 \begin{conj}
   \
 	\begin{itemize}
 	\item $\rho ^\perp(G) + \rho ^\perp(\bar{G}) \geq |V(G)| -2 $
 	\item $\rho ^\perp(G) + \rho ^\perp(\bar{G}) \leq |V(G)| +2 $
 	\end{itemize}
 \end{conj}
 {\em proof:} We based these conjectures in previous result, but \textcolor{red}{we have not demo}. we only know that in Cycles and Holes it's true and second sentence is false in $\sigma_{pg}$.
\section{Kernel of the dimension and dominant sets}
\begin{prop}
If $ \bar{G} $ is k-connected and $ V_c$ (cut vertex set) is independent set, then $ \rho ^ \perp (G) \le \rho ^ \perp (A) + \rho ^ \perp (B) + K $.
\end{prop}

\begin{prop}
If $ \bar {G} $ is k-connected and $V_c $(cut vertex set) is independent set, then $ \rho ^ \perp (G) \le \rho ^ \perp (A) + \rho ^ \perp (B) + \rho ^ \perp (V_c) $.
\end{prop}

{\em trivial}
$$\rho^\perp(A) + \rho^\perp(B) \le \rho^\perp(G)$$
\begin{conj} If $ \rho ^ \perp (V_c)> \rho ^ \perp (A) + \rho ^ \perp (b) $, then $$ \rho ^ \perp (G) = \rho ^ \perp (V_c) $$
\end{conj}
The Smollyn10 is a counterexample for greater than or equal to the conjecture.
\begin{conj}
It is true  $ \rho ^ \perp (G)> \rho ^ \perp (V_c) $ if it met the following:
\begin{itemize} \item $ \rho ^ \perp (V_c)> \rho ^ \perp (A) $ \item $ \rho ^ \perp (V_c)> \rho ^ \perp (B) $ \item $ \rho ^ \perp (V_c) \le \rho ^ \perp (A) + \rho ^ \perp (B) $
\end{itemize}
\end{conj}
\begin{conj} \textcolor {red} {FALSE: counterexample $ \bar {C_{11}} $}. \\If $ \rho ^ \perp (A) \ge \rho ^ \perp (V_c) $ and $ \rho ^ \perp (B) \ge \rho ^ \perp (V_c) $, then $$ \rho ^ \perp (G) = \rho ^ \perp (A) + \rho ^ \perp (B) $$
\end{conj}

\section{Relationship between different ORs}
Using the Lovasz work about orthogonal representation is trivial prove that if exists $\rho ^\perp(G)$ then exists $\rho(G)$, so $\rho^\perp(G) \geq \rho(G)$.
Also is trivial that $\rho_{lo}(G)$(not faithful) $\geq \rho(G)$ and $\rho_{pg}(G) \geq \rho (G)$. By simulations and computer algorithms we know that in general $\rho_{pg}(G) \geq \rho_{lo}(G)$ and moreover $\rho^\perp_{lo}(G) \geq \rho_{lo}(G)$. Then
$$\rho_{pg}(G) \geq \rho^\perp_{lo}(G) \geq \rho_{lo}(G) \geq \rho(G)$$
The most important parameters in our study are $\rho^\perp(G)$, $\rho^\perp_{pg}(G)$ and $\rho^\perp_{lo}(G)$.  There are a few known relationships between them based on their restrictions.
$$\rho^\perp_{pg} \geq \rho^\perp(G) \geq \rho(G)$$
$$\rho^\perp_{lo} \geq \rho^\perp(G) \geq \rho(G)$$
And the distance
$$max(\rho^\perp_{pg}(G)- \rho^\perp(G) ) = ||\rho^\perp_{pg}(G)- \rho^\perp(G)|| = |V(G)|-3$$


\begin{thebibliography}{99}

\bibitem{1}
Lovasz, L., Saks, M., Shrijver, A. Orthogonal representation and conectivity of graph, {\it Linear Algebra and its applications}, {\bf 4} (1987), no. 114/115, 439--454.

\bibitem{2}
Bell, John S. On the einstein-podolsky-ORsen paradox, {\it Physics}, {\bf 1} (1964), no. 3, 195--200.

\bibitem{3}
Cabello, Adán, Eirik Danielsen, Lars, López-Tarrida, Antonio J., Portillo, José R. Basic exclusivity graphs in quantum correlations, {\it Submmited in Physical Review Adwars}, {\bf } (2013).

\bibitem{4}Bell, John S. Bell, Rev.   Mod. Phys. {\bf 38} (1966) 447.

\bibitem{5} Kochen, S., Specker, E.P.  The problem of hidden variables in quantum mechanics, {\it Journal of Mathematics and Mechanics}, {\bf 17} (1967) 59-87.

\bibitem{6} Cabello, A., M. Estebaranz, José, Garcia-Alcaine, Guillermo. Bell-Kochen-Specker theorem: A proof with 18 vectors, {\it Physics letters} {\bf 212} (1996) 183-187 .

\bibitem{7} Peres, Asher. Quantum Theory: Concepts and Methods, {\it Kluwer Academic Publishers}, 2002.

\bibitem{8} Lovász, László. On the Shannon Capacity of a Graph, {\it IEEE Transactions on Infromation Theory}, {\bf 1} (1979) 25.

\bibitem{9} Alon, N. The Shannon Capacity of a Union, {\it Combinatorica} {\bf 18} (1998) 301-310.

\bibitem{10} Cabello, A., Severini, Simone, Winter, Andreas. (Non-)Contextuality of Physical Theories as an Axiom. arXiv:1010.2163v1. (2010).

\bibitem{11} Amselem, E., Eirik Danielsen, L.,  López-Tarrida, Antonio J., Portillo, Jose R., Bourennane, M., Cabello, A. Experimental Fully Contextual Correlations. arXiv:1111.3743v3 (2012).

\bibitem{12} Parsons, Pisansky. Vector representation of graphs, {\it Discrete Mathematics} {\bf 78} (1989) 143-154.

\bibitem{13} Sinajova, E. A note on vector representation of graphs, {\it Discrete Mathematics} {\bf 89} (1989) 315-317.

\bibitem{14} Alekseev, V.E., Lozin, V. On orthogonal representations of graphs, {\it Rutcor Research Report} (2000).

\bibitem{15} Neos Server Web Page. University of Wisconsin - Madison. http://www.neos-server.org/neos/solvers/index.html.

\bibitem{16} Solís, Alberto, Algoritmos para la resolución del problema de Representación Ortogonal, {\it Trabajo fin de máster en la Universidad de Sevilla}, Departamento de Matemática Aplicada I (2012).

\bibitem{17} Kennedy, J., Eberhart, R.C. Particle swarm optimization. {\it Proceedings of IEEE International Conference on Natural Networks, Perth, Australia} (1995) 1942-1948.

\bibitem{18} Lu, Haiyan, Chen, Weiqi. Self-adaptative velocity particle swarm optimization for solving constrained optimization problems. {\it J Glob Optim} {\bf 41} (2007) 427-445.

\bibitem{19} Trelea, Ioan C. The particle swarm optimization algorithm: convergence analysis and parameter selection. {\it Infromation processing letters} {\bf 85} (2003) 317-325.

\bibitem{code1}
 E. P. Specker,
 \href{http://onlinelibrary.wiley.com/doi/10.1111/j.1746-8361.1960.tb00422.x/abstract}{Dialectica \textbf{14}, 239 (1960);}
 \href{http://arxiv.org/abs/1103.4537}{arXiv:1103.4537.}

\bibitem{code2}
S.Kochen and E.P. Specker, in {\it Symposium on the Theory of Models}, edited by J. W. Addison, L. Henkin, and A. Tarski (North-Holland, Amsterdam, Holland, 1965), p. 177.

\bibitem{code3}
J.S Bell, \href{http://journals.aps.org/rmp/abstract/10.1103/RevModPhys.38.447}{Rev. Mod. Phys. {\bf 38}, 447 (1966).}

\bibitem{code4}
 S. Kochen and E. P. Specker,
 \href{http://www.iumj.indiana.edu/IUMJ/fulltext.php?year=1968&volume=17&artid=17004}{J. Math. Mech. \textbf{17}, 59 (1967).}

\bibitem{code5}
A. Stairs, \href{http://www.jstor.org/discover/10.2307/187557?uid=3737952&uid=2&uid=4&sid=21103727445967}{Phil. Sci. {\bf 50}, 578 (1983).}

\bibitem{code6}
R. Clifton, \href{http://scitation.aip.org/content/aapt/journal/ajp/62/5/10.1119/1.17551}{Am. J. Phys. {\bf 61}, 443 (1993)}; H. Bechmann Johansen, Am. J. Phys. {\bf 62}, 471 (1994); P.E. Vermaas, \href{http://scitation.aip.org/content/aapt/journal/ajp/62/7/10.1119/1.17488}{Am. J. Phys. {\bf 62}, 658 (1994).}

\bibitem{code7}
A. Cabello and G. García-Alcaine, \href{http://iopscience.iop.org/0305-4470/28/13/016/}{J. Phys. A \textbf{28}, 3719 (1995).}

\bibitem{code8}
A. Cabello and M. Terra Cunha, 
 \href{http://arxiv.org/pdf/1009.2330}{arXiv:1009.2330.}

\bibitem{code9}
A. Cabello and G. García-Alcaine, \href{http://iopscience.iop.org/0305-4470/29/5/016/}{J. Phys. A \textbf{29}, 1025 (1996).}

\bibitem{code10}
A. Cabello, \href{http://journals.aps.org/prl/abstract/10.1103/PhysRevLett.95.210401}{Phys. Rev. Lett. \textbf{101}, 210401 (2008).}
\href{http://xxx.lanl.gov/abs/quant-ph/0507259}{Cornell University Library, 0507259.}

\bibitem{code11}
P. Badziag, I. Bengtsson, A. Cabello, and I. Pitowsky, \href{http://journals.aps.org/prl/abstract/10.1103/PhysRevLett.103.050401}{Phys. Rev. Lett. \textbf{103}, 050401 (2009).}
\href{http://arxiv.org/abs/0809.0430}{arXiv:0809.0430.}

\bibitem{code12}
G. Kirchmair, F. Z\"{a}hringer, R. Gerritsma, M. Kleinmann, O. G\"{u}hne, A. Cabello, R. Blatt, and C. F. OORs, \href{http://www.nature.com/nature/journal/v460/n7254/full/nature08172.html}{Nature (London) \textbf{460}, 494 (2009).}
\href{http://arxiv.org/abs/0904.1655}{arXiv:0904.1655.}

\bibitem{code13}
E. Amselem, M. R{\aa}dmark, M. Bourennane, and A. Cabello, \href{http://journals.aps.org/prl/abstract/10.1103/PhysRevLett.103.160405}{Phys. Rev. Lett. \textbf{103}, 160405 (2009).}
\href{http://arxiv.org/abs/0907.4494}{arXiv:0907.4494.}

\bibitem{code14}
A. Cabello, \href{http://journals.aps.org/prl/abstract/10.1103/PhysRevLett.104.220401}{Phys. Rev. Lett. \textbf{104}, 220401 (2010).}
\href{http://arxiv.org/abs/0910.5507}{arXiv:0910.5507.}

\bibitem{code15}
B.D. MacKay, \texttt{nauty} \href{http://cs.anu.edu.au/~Brendan.McKay/nauty/nug.pdf}{User's Guide (Version 2.4)} (Departament of Computer Science, Australian National University, Canberra, Australia, 2007).

\bibitem{code16}
A. Peres, \emph{Quantum Theory: Concepts and Methods} (Kluwer, Dordrecht, 1993), p. 209.

\bibitem{code17}
A. Cabello, J. M. Estebaranz and G. García-Alcaine, \href{http://www.sciencedirect.com/science/article/pii/037596019600134X?via=ihub}{Phys. Rev. Lett. A \textbf{212}, 183 (1996), Eqs. (20)-(26).}

\bibitem{code18}
A. Cabello, \emph{Pruebas Algebraicas de Imposibilidad de Variables Ocultas en Mec{á}nica Cu{á}ntica}, Ph. D. Thesis, Universidad Complutense de Madrid, 1996, p. 201. It can be downloaded from \href{http://www.adancabello.com}{www.adancabello.com.}

\bibitem{code19}
A. Peres, \href{http://link.springer.com/article/10.1023\%2FA\%3A1026000614638}{Found. Phys. \textbf{33}, 1543 (2003).}
\href{http://arxiv.org/pdf/quant-ph/0207020}{arXiv:0207020.}

\bibitem{code20}
A. Cabello, J. R. Portillo, and G. Potel (unpublished).

\bibitem{code21}
A. Cabello, L.E. Danielsen, A.J. López-Tarrida, and J.R. Portillo, \emph{Basic exclusivity graphs in quantum correlations},  \href{http://dx.doi.org/10.1103/PhysRevA.88.032104}{Phys. Rev. A 88, 032104 ñ Published 9 September 2013.}
\href{http://arxiv.org/abs/1211.5825}{arXiv:1211.5825.}

\bibitem{41}
Leslie Hogben, \textit{Orthogonal representations, minimum rank, and graph complements}, Linear Algebra and its Applications, Volume 428, Issues 11–12, 1 June 2008, Pages 2560-2568, ISSN 0024-3795, \href{http://dx.doi.org/10.1016/j.laa.2007.12.004}{http://dx.doi.org/10.1016/j.laa.2007.12.004}.

\bibitem{42}
Chathand and Harary, 1968

\bibitem{43}
American Institute of Mathematics, Minimum rank graph catalog. (\href{http://aimath.org/pastworkshops/
matrixspectrum.html}{http://aimath.org/pastworkshops/
matrixspectrum.html})

\bibitem{44}
Jia–De Lin, Yue–Li Wang, Jou–Ming Chang and Hung–Chang Chan. \textit{On the Connectivity of a Graph and its Complement Graph}. \href{http://citeseerx.ist.psu.edu/viewdoc/download?doi=10.1.1.381.7661\&rep=rep1\&type=pdf}{http://citeseerx.ist.psu.edu/viewdoc/download?doi=10.1.1.381.7661\&rep=rep1\&type=pdf}

\bibitem{45}
Barioli, Francesco; Fallat, Shaun; and Hogben, Leslie (2005) "A variant on the graph parameters of Colin de Verdiere: Implications to the minimum rank of graphs," Electronic Journal of Linear Algebra: Vol. 13, Article 24. 
\href{http://dx.doi.org/10.13001/1081-3810.1170}{DOI: http://dx.doi.org/10.13001/1081-3810.1170}


\end{thebibliography}
\end{document}